\documentclass[11pt,a4paper]{amsart}

\title[Analytic continuation of residue currents]{Analytic continuation of residue currents}

\author{H\aa kan Samuelsson}
\thanks{Author supported by a Post Doctoral Fellowship from the Swedish Research Council.}
\address{Department of Mathematics, University of Wuppertal, Gaussstrasse 20, 42119 Wuppertal, Germany}
\email{hasam@math.chalmers.se}

\usepackage[english]{babel}
\usepackage{amsmath}
\usepackage{amsthm,amssymb,latexsym}
\usepackage{t1enc}
\usepackage{fancyheadings}
\usepackage{mathrsfs}
\usepackage{graphicx}

\newtheorem{proposition}{Proposition}
\newtheorem{theorem}[proposition]{Theorem}
\newtheorem{lemma}[proposition]{Lemma}

\theoremstyle{definition}

\newtheorem{remark}[proposition]{Remark}

\newcommand{\C}{\mathbb{C}}
\newcommand{\N}{\mathbb{N}}
\newcommand{\debar}{\bar{\partial}}
\newcommand{\D}{\mathscr{D}}
\newcommand{\R}{\mathbb{R}}

\newcommand{\Q}{\mathbb{Q}}

\newcommand{\CP}{\mathbb{CP}}

\newcommand{\real}{\mathfrak{R}\mathfrak{e}}

\def\newop#1{\expandafter\def\csname #1\endcsname{\mathop{\rm #1}\nolimits}}
\newop{span}

\begin{document}
\nocite{*}
\bibliographystyle{plain}

\begin{abstract}
Let $X$ be a complex manifold and $f\colon X\rightarrow \C^p$ a holomorphic mapping defining 
a complete intersection. We  
prove that the iterated Mellin transform of the residue
integral associated to $f$ has an analytic continuation to a neighborhood of the origin
in $\C^p$. 
\end{abstract}

\maketitle
\thispagestyle{empty}

\section{Introduction}
Let $X$ be a complex manifold of complex dimension $\textrm{dim}_{\C}X=n$ and let 
$f = (f_1,\ldots,f_{p+q}) \colon X \rightarrow \C^{p+q}$ be a holomorphic mapping
defining a complete intersection.
For a test form $\varphi \in \D_{n,n-p}(X)$ we let the residue integral of $f$ be the 
integral
\begin{equation*}
I_f^{\varphi}(\epsilon)=
\int_{T(\epsilon)}\frac{\varphi}{f_1\cdots f_{p+q}},
\end{equation*} 
where $T(\epsilon)$ is the tubular set 
$T(\epsilon)=\cap_1^p\{|f_j|^2=\epsilon_j\}\bigcap\cap_{p+1}^{p+q}\{|f_j|^2>\epsilon_j\}$.
If we let $\epsilon$ tend to zero along a path to the origin in the
first orthant such that $\epsilon_j/\epsilon_{j+1}^{k}\rightarrow 0$ for $j=1,\ldots,p+q-1$ and all
$k\in \N$, a so called ``admissible path'', 
then by fundamental results of Coleff-Herrera, \cite{CH}, and Passare, \cite{PCrelle},
the residue integral converges and the limit defines the action of a $(0,p)$-current
on the test form $\varphi$. We will refer to this current as the Coleff-Herrera-Passare current
and denote it suggestively by
\begin{equation*}
\debar \big[\frac{1}{f_1}\big]\wedge \cdots \wedge \debar \big[\frac{1}{f_p}\big]
\big[\frac{1}{f_{p+1}}\big]\cdots \big[\frac{1}{f_{p+q}}\big],
\end{equation*}
or sometimes $R^pP^q[1/f]$ for short. The current $R^p[1/f]$ is the classical Coleff-Herrera
residue current, which has proven to be a good notion of a multivariable residue of $f$, but
also the currents $R^pP^q[1/f]$, with $q\geq 1$, have turned out to be important to the theory.
In particular, if $q=1$ then $R^pP^q[1/f]$ is a $\debar$-potential to the Coleff-Herrera
residue current. 

The first question raised by Coleff and Herrera in the book \cite{CH} is whether
the residue integral $I_f^{\varphi}(\epsilon)$ has an unrestricted limit as $\epsilon$ tends to
zero. This question was answered in the negative by Passare and Tsikh in \cite{PTmotex}, where
they found two polynomials in $\C^2$, with the origin as the only common zero, such that 
the corresponding residue integral does not converge unrestrictedly; large classes of such examples 
were then found by Bj\"{o}rk. In this sense, the definition of the Coleff-Herrera-Passare current
is quite unstable. A different and, as we will see, more rigid approach is based on 
analytic continuation. Let $\lambda_1,\ldots,\lambda_{p+q}$ be complex parameters with 
$\real \, \lambda_j$ large. Then the integral
\begin{equation*}
\Gamma_f^{\varphi}(\lambda)=
\int_X \frac{\debar |f_1|^{2\lambda_1}\wedge \cdots \wedge \debar |f_p|^{2\lambda_p} |f_{p+1}|^{2\lambda_{p+1}}
\cdots |f_{p+q}|^{2\lambda_{p+q}}}{f_1\cdots f_{p+q}}\wedge \varphi,
\end{equation*}
makes sense and defines an analytic function of $\lambda$. This function is the iterated Mellin transform of the 
residue integral, i.e., 
\begin{equation*}
\pm \Gamma_f^{\varphi}(\lambda)=
\int_{[0,\infty)^{p+q}} I_f^{\varphi}(s) \, d(s_1^{\lambda_1})\wedge \cdots \wedge d(s_{p+q}^{\lambda_{p+q}}).
\end{equation*}
It is known that $\Gamma_f^{\varphi}(\lambda)$ has a meromorphic continuation to all of $\C^{p+q}$ and that 
its only possible poles in a neighborhood of the half space $\cap_1^{p+q}\{\real \,\lambda_j\geq 0\}$
are along hyperplanes of the form $\sum a_j\lambda_j=0$, $a_j\in \Q_+$. Moreover, by results of Yger,
the restriction of $\Gamma_f^{\varphi}(\lambda)$ to any complex line of the form 
$\{\lambda=(t_1z,\ldots,t_pz);\, z\in \C\}$, $t_j\in \R_+$, is analytic at the origin and the value 
there equals the action of the Coleff-Herrera-Passare current on $\varphi$. In the case of codimension two,
i.e., when $f=(f_1,f_2)$, it is also known
that the corresponding $\Gamma$-functions are analytic 
in some neighborhood of $\cap_1^{2}\{\real \,\lambda_j\geq 0\}$. Yger has posed the question
whether this generalizes to arbitrary codimensions. 
The purpose of this paper is to prove the following theorem, which answers Yger's question
in the affirmative.

\begin{theorem}\label{analytisk}
Let $X$ be a complex manifold of complex dimension $n$ and 
let $f=(f_1,\ldots,f_{p+q}) \colon X \rightarrow \C^{p+q}$ be a holomorphic mapping defining a 
complete intersection. If $N$ is a positive integer and $\varphi\in \D_{n,n-p}(X)$ is 
a test form on $X$ then the integral
\begin{equation}\label{intro0}
\Gamma_{f^N}^{\varphi}(\lambda)=
\int_X \frac{\debar |f_1|^{2\lambda_1}\wedge \cdots \wedge \debar |f_p|^{2\lambda_p} |f_{p+1}|^{2\lambda_{p+1}}
\cdots |f_{p+q}|^{2\lambda_{p+q}}}{f_1^N\cdots f_{p+q}^N}\wedge \varphi,
\end{equation} 
is analytic in a half
space $\{\lambda \in \C^{p+q};\, \real\, \lambda_j>-\epsilon, \, 1\leq j \leq p+q\}$ 
for some $\epsilon\in \Q_+$ independent of $N$. 
\end{theorem}
We remark that for non-complete intersections, the $\Gamma$-function
still is meromorphic in $\C^{p+q}$ but will in general have poles along hyperplanes through the origin.
Our proof of Theorem \ref{analytisk} uses Hironaka's theorem
on resolutions of singularities, \cite{H}, to reduce to the case when $\{f_1\cdots f_{p+q}=0\}$ has normal
crossings, i.e., in local coordinates on a blow-up manifold lying above $X$, 
the pull-back, $\hat{f}_j$, of the $f_j$ are monomials times invertible holomorphic functions.
(For our proof it is actually enough to use the weaker version of Hironaka's theorem where the projection from the
blow-up manifold to $X$ is allowed to be ``finite to one'' outside the exceptional divisor.)
In general, the $\hat{f}_j$ do not define a complete intersection on the blow-up manifold but 
the information that the $f_j$ do on the base manifold is coded in the pull-back, $\hat{\varphi}$, of the 
test form $\varphi$. We are able to recover this information using a Whitney-type division lemma for
(anti-)holomorphic forms. It is also worth noticing that for $p=1$, the problem of analytic continuation
of $\Gamma_{f^N}^{\varphi}(\lambda)$ is of local nature on the blow-up manifold, i.e., 
it suffices to consider one chart on
the blow-up at the time. This is not the case if $p\geq 2$ and $q\geq 1$, all charts on the blow-up 
have to be considered simultaneously. We give a simple example showing this in Section \ref{ex}. 
In \cite{HS} we were able to overcome this problem in the special case when $p=2$ and $q=1$ by quite involved
integrations by parts on the blow-up manifold. Very rewarding discussions with Jan-Erik Bj\"{o}rk have 
resulted in a much simpler and more transparent argument based on induction over $p$. 

We continue and give a short historical account of analytic continuation of residue currents.
The case $p=0$ and $q=1$ is the most studied one and the analytic continuation was 
in this case proved by Atiyah in \cite{At} using Hironaka's theorem; see also \cite{BernGelf}.
The main point was to get a
multiplicative inverse of $f$ in the space of currents, and indeed,
the value at $\lambda=0$ gives a current, $U$, such that $fU=1$ in the sense of currents.
At the same time, Dolbeault and Herrera-Lieberman proved, also using Hironaka's theorem, that the principal
value current of $1/f$, defined by
\begin{equation*}
\D_{n,n}(X)\ni \varphi \mapsto \lim_{\epsilon\rightarrow 0}
\int_{|f|^2>\epsilon}\varphi/f,
\end{equation*}
and denoted $[1/f]$, exists, cf.\ \cite{Dol} and \cite{HL}. 
It is elementary to see that this current 
coincides with the current defined by Atiyah if $f$ is a monomial and for general $f$ it then follows 
from Hironaka's theorem. A perhaps more conceptional explanation for this equality is that 
the two definitions are linked via the Mellin transform; recall from above that 
$\int |f|^{2\lambda}\varphi/f$ is the Mellin transform of 
$\epsilon \mapsto \int_{|f|^2>\epsilon}\varphi/f$. The poles of the
current valued function $\lambda \mapsto |f|^{2\lambda}/f$
are closely related to the roots of the Bernstein-Sato polynomial, $b(\lambda)$, associated to $f$. 
By Bernstein-Sato theory, see, e.g., \cite{JEB}, $f$ satisfies some functional equation
\begin{equation*}
b(\lambda)\bar{f}^{\lambda}=\sum_j \lambda^j Q_j(\bar{f}^{\lambda+1}),
\end{equation*}
where $Q_j$ are anti-holomorphic differential operators. By iterating $m$ times and multiplying
with $f^{\lambda}/f^N$ it follows that
\begin{equation*}
b(\lambda+m)\cdots b(\lambda)|f|^{2\lambda}/f^{N}=\sum_j \lambda^j R_j(\bar{f}^{m}|f|^{2\lambda}/f^N)
\end{equation*} 
for some anti-holomorphic differential operators $R_j$. If 
$\varphi\in \D_{n,n}(X)$ and $R_j^*$ is the adjoint operator of $R_j$ it thus follows that
\begin{equation}\label{intro1}
\int_X |f|^{2\lambda}\frac{\varphi}{f^N}=b(\lambda+m)^{-1}\cdots b(\lambda)^{-1}
\sum_j \lambda^j \int_X |f|^{2\lambda}\frac{\bar{f}^m}{f^N}R_j^*(\varphi).
\end{equation}  
Now, from Kashiwara's result, \cite{Kash}, we know that $b(\lambda)$ has all of its roots contained 
in the set of negative rational numbers. 
Hence, we can read off from \eqref{intro1} that the current valued function
$\lambda \mapsto |f|^{2\lambda}/f^N$ has a meromorphic continuation to all of $\C$ and that its poles are 
contained in arithmetic progressions of the form $\{-s-\N\}$ with $s\in \Q_+$. In particular, 
$\lambda \mapsto \int |f|^{2\lambda}\varphi/f^N$ is holomorphic in some half space $\real \, \lambda >-\epsilon$
for some $\epsilon \in \Q_+$ independent of $N$. A detailed study of the poles that actually appear was done by 
Barlet in \cite{Bar}.
Consider now instead the current valued function $\lambda \mapsto \debar |f|^{2\lambda}/f$. It is the $\debar$-image
of $\lambda \mapsto |f|^{2\lambda}/f$ and has thus also a meromorphic continuation to all of $\C$ with poles contained
in arithmetic progressions of the form $\{-s-\N\}$. The value at $\lambda=0$ is now the residue current 
$\debar[1/f]$, i.e., the $\debar$-image of $[1/f]$. 
The case when $f$ is one function, i.e., $p+q=1$, is thus well understood.
When $p+q>1$, the picture is not that coherent. We have seen that $\Gamma_f^{\varphi}(\lambda)$
is the iterated Mellin transform of the residue integral
but from the examples by Passare-Tsikh and Bj\"{o}rk mentioned above, we know that $\Gamma_f^{\varphi}(\lambda)$
is {\em not} the Mellin transform of a continuous function in general.
A multivariable Bernstein-Sato approach
has been considered, but by results of Sabbah, \cite{Sabbah}, 
the zero set of the multivariable Bernstein-Sato polynomial  
will in general intersect $\cap_j \{\real \, \lambda_j\geq 0\}$, and so 
this method cannot be used to prove our result. On the other hand, it shows that 
$\Gamma_f^{\varphi}(\lambda)$ has a meromorphic continuation to all of $\C^{p+q}$.
More direct approaches have been considered by, e.g., Berenstein, Gay, Passare, Tsikh, and
Yger and the case $q=0$ has got the most attention. 
For instance, a direct proof of the meromorphic continuation of $\Gamma_f^{\varphi}(\lambda)$
to all of $\C^{p}$ can be found in \cite{PTcanada}. Also,
as mentioned above, it is proved in \cite{Y} that
the restriction of $\lambda \mapsto \Gamma_f^{\varphi}(\lambda)$
to any complex line of the form $\{\lambda=(t_1z,\ldots,t_pz);\, z\in \C\}$, where $t_j\in \R_+$,
is analytic at the origin and that the value there equals the Coleff-Herrera residue current. 
The first analyticity result in several variables was obtained by Berenstein
and Yger. They proved that if $p+q=2$, then $\Gamma_f^{\varphi}(\lambda)$ is in fact analytic in a 
half space in $\C^2$ containing the origin, see, e.g., \cite{BGVY} and \cite{PTcanada} for proofs. 
In view of these positive results it has been believed that this holds in general but to our knowledge,
no complete proof has appeared.

This paper is organized as follows.  
Section \ref{sketch} contains an outline
of the proof in the case $p=2$ and $q=1$. This is to show the essential steps without confronting 
technical and notational difficulties.
In Section \ref{ex} we give a simple example showing that global effects 
on the blow-up manifold have to be taken into account when $p\geq 2$ and $q\geq 1$.
The detailed proof of Theorem \ref{analytisk} is contained in Section \ref{bevis}.

\smallskip

{\bf Acknowledgment:} I am grateful to Jan-Erik Bj\"{o}rk for his help and support
during the preparation of this paper. His insightful comments and suggestions have substantially
improved and simplified many arguments as well as the presentation.

\section{The main elements of the proof}\label{sketch}
In this section we illustrate the main new ideas in our proof by considering the case when 
$p=2$ and $q=1$. Let $f_1$, $f_2$, and $f_3$ be holomorphic functions in $\C^3$ (for simplicity)
and assume that the origin is the only common zero. Using the techniques of, e.g., \cite{BGVY} or 
\cite{PTcanada} it is not hard to prove that the current valued function
$\lambda \mapsto (f_1^{-1}\debar |f_1|^{2\lambda_1}) f_2^{-1}|f_2|^{2\lambda_2}f_3^{-1}
|f_3|^{2\lambda_3}$ can be analytically continued to a neighborhood of
the origin; see also Proposition \ref{polprop} below. 
Assume now that we can prove that there is a polynomial $P_{12}(\lambda_1,\lambda_2)$,
which is a product of linear factors $a\lambda_1+b\lambda_2$ in $\lambda_1$ and $\lambda_2$,
such that the current valued function
\begin{equation}\label{sketch0}
\lambda \mapsto P_{12}(\lambda_1,\lambda_2)
\frac{\debar |f_1|^{2\lambda_1} \wedge \debar |f_2|^{2\lambda_2} |f_3|^{2\lambda_3}}{f_1 f_2 f_3}
\end{equation}
can be analytically continued to a neighborhood of the origin. That is, we assume for the moment 
that the only possible poles
(close to the origin) of the meromorphic current valued function $(f_1^{-1}\debar |f_1|^{2\lambda_1}) \wedge 
(f_2^{-1}\debar |f_2|^{2\lambda_2})f_3^{-1} |f_3|^{2\lambda_3}$ are along hyperplanes of the form
$a\lambda_1+b\lambda_2=0$. Consider the equality of currents
\begin{eqnarray}\label{sketch1}
\debar 
\frac{\debar |f_1|^{2\lambda_1} |f_2|^{2\lambda_2} |f_3|^{2\lambda_3}}{f_1 f_2 f_3} &=&
-\frac{\debar |f_1|^{2\lambda_1} \wedge \debar |f_2|^{2\lambda_2} |f_3|^{2\lambda_3}}{f_1 f_2 f_3} \\
& &
-\frac{\debar |f_1|^{2\lambda_1} |f_2|^{2\lambda_2} \wedge \debar |f_3|^{2\lambda_3}}{f_1 f_2 f_3}, \nonumber
\end{eqnarray}
which holds for $\real \, \lambda_1, \real\, \lambda_2, \real\,\lambda_3 >> 1$. We know that the left 
hand side can be analytically continued to a neighborhood of the origin and we have assumed that we
can prove that the first term on the right hand side only has poles (close to the origin) along hyperplanes
$a\lambda_1+b\lambda_2=0$. The last term on the right hand side therefore also has only such poles. But, by
permuting the indices, we can, assumingly, prove that the last term on the right hand side 
only has poles along hyperplanes of the form
$a'\lambda_1+b'\lambda_3=0$. We can thus conclude that the only possible pole that the last term on the 
right hand side can have is along $\lambda_1=0$. On the other hand,
if we switch the indices
$1$ and $3$ in \eqref{sketch1} we similarly get that the last term in \eqref{sketch1} only has poles along 
hyperplanes $a''\lambda_2+b''\lambda_3=0$. (The last term is unaffected by the switch modulo a sign.)
Its possible pole along $\lambda_1=0$ is thus not present.
In conclusion, the last term in \eqref{sketch1} has an analytic continuation to a neighborhood of the 
origin if we can prove the existence of a polynomial $P_{12}(\lambda_1,\lambda_2)$ such that 
\eqref{sketch0} can be analytically continued to a neighborhood of the origin. 
To do this, we use Hironaka's theorem to compute $\Gamma_f^{\varphi}(\lambda)$ on a blow-up manifold. 
More precisely, for some neighborhood 
$U$ of an arbitrary point in $\C^3$ one can find a blow-up manifold $\mathcal{U}$, lying 
properly above the base $U$, such that the preimage, $\mathcal{Z}$, of $Z=\{f_1f_2f_3=0\}$ 
has normal crossings and $\mathcal{U}\setminus \mathcal{Z}$ is biholomorphic with
$U\setminus Z$. By a partition of unity we may assume that $\varphi$ has support in such a $U$ and we 
pull our integral $\Gamma_f^{\varphi}(\lambda)$ back to $\mathcal{U}$. In local charts
on $\mathcal{U}$, we then have that the pullback, $\hat{f}_j$, of the $f_j$ are monomials, 
$x^{\alpha(j)}$, times invertible holomorphic functions. Let us consider a generic chart where
the multiindices $\alpha(1)$, $\alpha(2)$, and $\alpha(3)$ are linearly independent. It is then possible to define
new coordinates, still denoted $x$, such that the invertible holomorphic functions are $\equiv 1$;
see, e.g., \cite{PCrelle}. We note that, in general, there are also so called charts of resonance where
one cannot choose coordinates so that the invertible functions are $\equiv 1$. These charts are responsible
for the discontinuity of the residue integral, $I_f^{\varphi}(\epsilon)$, but do not cause any problems 
in our situation. 
We shall thus consider the integral
\begin{equation}\label{sketch2}
\int \frac{\debar |x^{\alpha(1)}|^{2\lambda_1}\wedge \debar |x^{\alpha(2)}|^{2\lambda_2}
|x^{\alpha(3)}|^{2\lambda_3}}{x^{\alpha(1)}x^{\alpha(2)}x^{\alpha(3)}}\wedge \rho \hat{\varphi},
\end{equation}   
where $\rho$ is some cut-off function on $\mathcal{U}$. In $\C^3$, we have 
$\varphi(z)=\sum_1^3\varphi_{j}(z)dz\wedge d\bar{z}_j$, and so, by linearity 
we may assume that $\varphi$ has a decomposition $\varphi=\phi\wedge \bar{\psi}$, where 
$\phi\in \D_{3,0}(\C^3)$ and $\bar{\psi}$ is the conjugate of a holomorphic $1$-form. The pullback
$\hat{\varphi}=\hat{\phi}\wedge \bar{\hat{\psi}}$ thus also has such a decomposition.
For simplicity we assume that $\hat{\psi}=h(x)dx_3$ for a holomorphic function $h$ on $\mathcal{U}$. 
Then \eqref{sketch2} equals
\begin{equation}\label{sketch3}
\lambda_1\lambda_2
\int \frac{|x^{\alpha(1)}|^{2\lambda_1} |x^{\alpha(2)}|^{2\lambda_2}
|x^{\alpha(3)}|^{2\lambda_3}}{x^{\alpha(1)}x^{\alpha(2)}x^{\alpha(3)}}
\frac{d\bar{x}^{\alpha(1)} \wedge d\bar{x}^{\alpha(2)}}{
\bar{x}^{\alpha(1)} \bar{x}^{\alpha(2)}}
\wedge \rho \hat{\varphi} \hspace{2cm}
\end{equation}
\begin{equation*}
\hspace{2cm} =
\lambda_1\lambda_2
\int \frac{|x^{\alpha(1)}|^{2\lambda_1} |x^{\alpha(2)}|^{2\lambda_2}
|x^{\alpha(3)}|^{2\lambda_3}}{x^{\alpha(1)}x^{\alpha(2)}x^{\alpha(3)}}
A_{12}\frac{d\bar{x}_1\wedge d\bar{x}_2}{\bar{x}_1\bar{x}_2}
\wedge \rho \hat{\varphi},
\end{equation*} 
where $A_{12}=\alpha(1)_1\alpha(2)_2-\alpha(1)_2\alpha(2)_1$. We may assume that $A_{12}>0$ and so,
in particular, $\alpha(1)_1>0$ and $\alpha(2)_2>0$. To avoid having to consider so many cases we also
assume that $\alpha(1)_2=\alpha(2)_1=0$. The cases when this is not fulfilled do not cause any 
additional difficulties and can be treated similarly. Three cases can occur:
\begin{itemize}
\item[i)] Neither $x_1$ nor $x_2$ divides $x^{\alpha(3)}$.
\item[ii)] Precisely one of $x_1$ and $x_2$ divides $x^{\alpha(3)}$.
\item[iii)] Both $x_1$ and $x_2$ divide $x^{\alpha(3)}$.
\end{itemize}
Consider the case ii) and assume that $x_1$ divides $x^{\alpha(3)}$. The variety
$V=\{f_1=f_3=0\}$ in $\C^3$ has codimension $2$ since $f_1$, $f_2$, and $f_3$ define a complete
intersection, and so the holomorphic $2$-form 
$df_2\wedge \psi$ has a vanishing pullback to $V$. Since $x_1$ divides both 
$\hat{f}_1=x^{\alpha(1)}$ and $\hat{f}_3=x^{\alpha(3)}$ we see that 
$d\hat{f}_2\wedge \hat{\psi}=dx^{\alpha(2)}\wedge \hat{\psi}$ must have a vanishing pullback to 
$\{x_1=0\} \subseteq \{\hat{f}_1=\hat{f}_3=0\}$. But $x_1$ does not divide
$x^{\alpha(2)}$ and hence, $dx^{\alpha(2)}\wedge \hat{\psi}|_{x_1=0}=0$ in $\mathcal{U}$, where 
$\hat{\psi}|_{x_1=0}$ means the pullback of $\hat{\psi}$ to $\{x_1=0\}$ extended constantly
to $\mathcal{U}$. This implies that we may replace $\hat{\varphi}=\hat{\phi}\wedge \bar{\hat{\psi}}$
in \eqref{sketch3} by $\hat{\phi}\wedge (\bar{\hat{\psi}}-\bar{\hat{\psi}}|_{x_1=0})$ without 
affecting the integral. Now, $x_1$ divides $\hat{\psi}-\hat{\psi}|_{x_1=0}$ so we may in fact assume 
that $\hat{\varphi}$ in \eqref{sketch3} is divisible by $\bar{x}_1$, or formulated differently,
that $(d\bar{x}_1/\bar{x}_1)\wedge \hat{\varphi}$ is a smooth form.
If we instead consider the case iii), similar degree arguments give that
$dx^{\alpha(1)}\wedge \hat{\psi}|_{x_2=0}= dx^{\alpha(2)}\wedge \hat{\psi}|_{x_1=0}=0$ in
$\mathcal{U}$. We may then replace $\hat{\varphi}$ in \eqref{sketch3} by 
\begin{equation}\label{sketch4}
\hat{\phi}\wedge (\bar{\hat{\psi}}-\bar{\hat{\psi}}|_{x_1=0}-
\bar{\hat{\psi}}|_{x_2=0}+\bar{\hat{\psi}}|_{x_1=x_2=0})
\end{equation}
without affecting the integral. But \eqref{sketch4} is divisible by $\bar{x}_1\bar{x}_2$,
and so, in this case, we may assume that 
$(d\bar{x}_1/\bar{x}_1)\wedge (d\bar{x}_2/\bar{x}_2)\wedge \hat{\varphi}$ is a smooth form.
It is now easy to see that \eqref{sketch3} has a meromorphic continuation and that its possible poles
close to the origin are along hyperplanes $a\lambda_1+b\lambda_2=0$. In case i), we write 
\eqref{sketch3} as
\begin{equation*}
A_{12} \frac{\lambda_1\lambda_2}{\mu_1\mu_2}
\int \frac{\debar |x_1|^{2\mu_1}\wedge \debar |x_2|^{2\mu_2} |x_3|^{2\mu_3}}{x_1^{\alpha_1}
x_2^{\alpha_2}x_3^{\alpha_3}}\wedge \rho \hat{\varphi},
\end{equation*}
where $\mu_j=\sum_{i=1}^3\lambda_i\alpha(i)_j$ and $\alpha_j=\sum_{i=1}^3\alpha(i)_j$.
It is an easy one-variable problem to see that this integral
(without the coefficient) has an analytic continuation to a neighborhood of the origin; cf., e.g.,
Lemma 2.1 in \cite{matsa}. 
Since neither $x_1$ nor $x_2$ divides $x^{\alpha(3)}$, i.e., $\alpha(3)_1=\alpha(3)_2=0$, we have 
$\mu_1=\alpha(1)_1\lambda_1+\alpha(2)_1\lambda_2$ and $\mu_2=\alpha(1)_2\lambda_1+\alpha(2)_2\lambda_2$
and it follows that \eqref{sketch3} only has poles of the allowed type in the case i). In the case 
ii) (with $x_1$ dividing $x^{\alpha(3)}$) we write \eqref{sketch3} as
\begin{equation*}
-A_{12}\frac{\lambda_1\lambda_2}{\mu_2}
\int \frac{ |x_1|^{2\mu_1}\debar |x_2|^{2\mu_2} |x_3|^{2\mu_3}}{x_1^{\alpha_1}
x_2^{\alpha_2}x_3^{\alpha_3}}\wedge \frac{d\bar{x}_1}{\bar{x}_1}\wedge \rho \hat{\varphi}.
\end{equation*}
Now $\mu_2=\alpha(1)_2\lambda_1+\alpha(2)_2\lambda_2$, since $x_2$ does not divide $x^{\alpha(3)}$,
and from our considerations above we may assume that $(d\bar{x}_1/d\bar{x}_1)\wedge \hat{\varphi}$ 
is smooth. It follows that 
\eqref{sketch3} only has the allowed type of poles in the case ii) as well. The case iii) is easier;
then we may assume that 
$(d\bar{x}_1/\bar{x}_1)\wedge (d\bar{x}_2/\bar{x}_2)\wedge \hat{\varphi}$ is a smooth form and 
\eqref{sketch3} is in this case even analytic at the origin. 

\begin{remark}
As mentioned in the introduction, and shown in the next section, it is necessary to take 
global effects on the blow-up manifold into account when proving analyticity of \eqref{intro0}
beyond the origin. However, as indicated by the above argument, the problem of showing that
\eqref{intro0} only has poles along hyperplanes of the form 
$\sum_1^pa_j\lambda_j=0$ is of a local nature on the blow-up; cf.\ \cite{HS}.
\end{remark}

\section{An example}\label{ex}
We present a simple example showing that global effects on the blow-up manifold have to be 
taken into account when proving our result for $p\geq 2$ and $q\geq 1$. Consider the integral
\begin{equation}\label{eq0}
\int \frac{|x_1|^{2\lambda_1}\debar |x_2|^{2\lambda_2}\wedge \debar
|x_3|^{2\lambda_3}}{x_1x_2x_3}\wedge \varphi(x) dx\wedge d\bar{x}_1
\end{equation}
in $\C^3$, where $\varphi$ is a function defined as follows. Let 
$\phi$, $\varphi_2$ and $\varphi_3$ be smooth functions on 
$\C$ with support close to the origin but non-vanishing there, and put
$\varphi_1=\partial\phi/\partial\bar{z}$. We define $\varphi(x)$ to be the function
$\varphi_1(x_1)\varphi_2(x_2)\varphi_3(x_3)$ in $\C^3$. 
Note that \eqref{eq0} equals 
\begin{equation*}
\int \frac{|x_1|^{2\lambda_1} |x_2|^{2\lambda_2}|x_3|^{2\lambda_3}}{x_1x_2x_3} 
\varphi_1 \frac{\partial \varphi_2}{\partial \bar{x}_2}
\frac{\partial \varphi_3}{\partial \bar{x}_3}dx\wedge d\bar{x}
\end{equation*}
after two integrations by parts, from which we see that \eqref{eq0} is analytic at $\lambda=0$.
Now we blow up $\C^3$ along the $x_1$-axis and look at the pullback of \eqref{eq0} to this 
manifold. Let $\pi\colon \C\times \mathcal{B}_0\C^2\rightarrow \C^3$ be the blow-up map.
In the natural coordinates $z$ and $\zeta$ on $\C\times \mathcal{B}_0\C^2$ it then looks like
\begin{equation*}
\pi(z_1,z_2,z_3)=(z_1,z_2,z_2z_3),
\end{equation*}
\begin{equation*}
\pi(\zeta_1,\zeta_2,\zeta_3)=(\zeta_1,\zeta_2\zeta_3, \zeta_2).
\end{equation*}
Since $\varphi$ has support close to the origin, $\pi^*\varphi$ has support close to
$\pi^{-1}(0)=\{z_1=z_2=0\}\cup\{\zeta_1=\zeta_2=0\}\cong \CP^1$. Note that $z_3$ and 
$\zeta_3$ are natural coordinates on this $\CP^1$ and choose a partition of unity,
$\{\rho_1,\rho_2\}$ on $\mbox{supp} (\pi^*\varphi)$ such that $\mbox{supp}
(\rho_1)\subset \{|z_3|<2\}$ and $\mbox{supp} (\rho_2)\subset \{|\zeta_3|<2\}$.
The pullback of \eqref{eq0} under $\pi$ now equals
\begin{equation*}
\int \frac{|z_1|^{2\lambda_1}\debar |z_2|^{2\lambda_2}\wedge \debar |z_2z_3|^{2\lambda_3}}{z_1z_2^2z_3}
\wedge \rho_1(z)\varphi_1(z_1)\varphi_2(z_2)\varphi_3(z_2z_3)z_2dz\wedge d\bar{z}_1
\end{equation*}
\begin{equation*}
-\int \frac{|\zeta_1|^{2\lambda_1}\debar |\zeta_2\zeta_3|^{2\lambda_2}\wedge
\debar |\zeta_2|^{2\lambda_3}}{\zeta_1\zeta_2^2\zeta_3}\wedge
\rho_2(\zeta)\varphi_1(\zeta_1)\varphi_2(\zeta_2\zeta_3)\varphi_3(\zeta_2)\zeta_2
d\zeta\wedge d\bar{\zeta}_1.
\end{equation*}
We know that this sum (difference) is analytic at $\lambda=0$ but we will see that none of the terms 
are. Consider the first term. It is easily verified that it can be written as
\begin{equation*}
\frac{\lambda_2}{\lambda_2+\lambda_3}
\int \frac{|z_1|^{2\lambda_1}\debar |z_2|^{2(\lambda_2+\lambda_3)}\wedge 
\debar |z_3|^{2\lambda_3}}{z_1z_2z_3}
\rho_1(z)\varphi_1(z_1)\varphi_2(z_2)\varphi_3(z_2z_3)dz\wedge d\bar{z}_1.
\end{equation*}
We denote this integral, with the coefficient $\lambda_2/(\lambda_2+\lambda_3)$ removed,
by $I(\lambda)$. After two integrations by parts one sees that $I(\lambda)$ is analytic at the
origin, and so $\lambda_2I(\lambda)/(\lambda_2+\lambda_3)$ is analytic at the origin if and 
only if $I(\lambda)$ vanishes on the hyperplane $\lambda_2+\lambda_3=0$. In particular
we must have that $I(0)=0$. But $I(0)$ can be computed using Cauchy's formula, and one obtains
$I(0)=-(2\pi i)^3\phi(0)\varphi_2(0)\varphi_3(0)\neq 0$. 

\begin{remark}
This example could be a little confusing. The variable $z_1$ just appears as a ``dummy variable'' 
in the computations above, to which nothing interesting happens. This indicates that global effects appear
already in the case $p=2$ and $q=0$. It is in fact so,
but in this case the analyticity follows, simply by applying to $\debar$-exact test forms,
if we can prove analyticity of 
\begin{equation*}
(\lambda_1,\lambda_2)\mapsto \int 
\frac{|f_1|^{2\lambda_1}\debar |f_2|^{2\lambda_2}}{f_1\cdot f_2}\wedge \varphi, \,\,\,\,
\varphi \in \D_{n,n-1}(X).
\end{equation*}
This can can actually be done using only local arguments, see, e.g., Proposition \ref{polprop} below.
The case $p=2$, $q=0$ can therefore be {\em reduced} to a case where only local arguments are needed,
but for $p\geq 2$ and $q\geq 1$ this is in general not possible. 
\end{remark}

\section{The proof}\label{bevis}
We give here the detailed proof of Theorem \ref{analytisk}. We begin with the following proposition,
whose proof relies on the Whitney-type division lemma (Lemma \ref{divlemma}) below.

\begin{proposition}\label{polprop}
Let $f=(f_1,\ldots,f_p)\colon X\rightarrow \C^p$ be a holomorphic mapping defining a complete 
intersection and let $g_1,\ldots,g_q$ be holomorphic functions on $X$ such that 
$(f_1,\ldots,f_p,g_j)$ defines a complete intersection for each $j=1,\ldots,q$. Let also 
$\varphi\in \D_{n,n-p}(X)$ be a test form and $N$ a positive integer. Then, for some $\epsilon\in \Q_+$
independent of $N$, the function
\begin{equation*}
\Gamma_{f,g}^{\varphi}(\lambda)=
\int_X\frac{\debar |f_1|^{2\lambda_1}\wedge \cdots \wedge \debar |f_p|^{2\lambda_p}
|g_1|^{2\lambda_{p+1}}\cdots |g_q|^{2\lambda_{p+q}}}{f_1^N\cdots f_p^N g_1^N \cdots g_q^N}\wedge \varphi,
\end{equation*}
originally defined when all $\real\, \lambda_j$ are large, has a meromorphic continuation to 
all of $\C^{p+q}$ and its only possible poles in the half space 
$H=\{\lambda\in \C^{p+q};\, \real\,\lambda_j >-\epsilon, \, 1\leq j \leq p+q\}$
are along hyperplanes of the form $\sum_1^pa_j\lambda_j=0$, where $a_j\in \N$
and at least two of the $a_j$ are non-zero. In particular, if $p=1$ then $\Gamma_{f,g}^{\varphi}(\lambda)$
is analytic in $H$.
\end{proposition}
\begin{proof}
It is well known that $\Gamma_{f,g}^{\varphi}(\lambda)$ has a meromorphic continuation to all of 
$\C^{p+q}$ so we only check that its possible poles in $H$ are of the prescribed form. 
We will compute $\Gamma_{f,g}^{\varphi}(\lambda)$ by pulling the integral back to a blow-up
manifold, $\mathcal{X}$, given by Hironaka's theorem, where the variety 
$\{\hat{f}_1\cdots \hat{f}_p\cdot \hat{g}_1
\cdots \hat{g}_q=0\}$ has normal crossings; cf.\ Section \ref{sketch}. 
(The hat, $\hat{\,\,\,}$, means pullback to the blow-up.)
We can thus write
\begin{equation*}
\Gamma_{f,g}^{\varphi}(\lambda)=
\sum_{\rho}\int_{\mathcal{X}}
\frac{\debar |\hat{f}_1|^{2\lambda_1}\wedge \cdots \wedge \debar |\hat{f}_p|^{2\lambda_p}
|\hat{g}_1|^{2\lambda_{p+1}}\cdots |\hat{g}_q|^{2\lambda_{p+q}}}{
\hat{f}_1^N\cdots \hat{f}_p^N \hat{g}_1^N \cdots \hat{g}_q^N}\wedge \rho \hat{\varphi},
\end{equation*}
where $\{\rho\}$ is a partition of unity of $\textrm{supp}(\hat{\varphi})$ and each $\rho$ has support 
in a coordinate chart where $\hat{f}_i$ and $\hat{g}_j$ are monomials times invertible holomorphic
functions. Let us consider a chart with holomorphic coordinates $x$ in which
$\hat{f}_1=u_1x^{\alpha(1)},\ldots,\hat{f}_p=u_px^{\alpha(p)}$ and 
$\hat{g}_1=v_1x^{\beta(1)},\ldots, \hat{g}_q=v_qx^{\beta(q)}$, 
where the $u_i$ and the $v_j$ are invertible and holomorphic. Denote by $m$ the number of 
vectors in a maximal linearly independent subset of $\{\alpha(1),\ldots,\alpha(p)\}$ and assume 
for simplicity that $\alpha(1),\ldots,\alpha(m)$ are linearly independent. It is then possible 
to define new coordinates, still denoted by $x$, such that $u_1=\cdots =u_m\equiv 1$ in the new 
coordinates; see, e.g., \cite{PCrelle} page 46. Now, for each $j=m+1,\ldots,p$, $\alpha(j)$
is a linear combination of $\alpha(1),\ldots,\alpha(m)$ and it follows from exterior algebra
that $dx^{\alpha(j)}\wedge dx^{\alpha(1)}\wedge \cdots \wedge dx^{\alpha(m)}=0$. In the 
$x$-chart, the term we are looking at can therefore be written
\begin{equation}\label{prop1}
\int
\frac{\debar |x^{\alpha(1)}|^{2\lambda_1}\wedge \cdots \wedge \debar |x^{\alpha(m)}|^{2\lambda_m}
|x^{\lambda\gamma}|^{2}}{x^{N\alpha}x^{N\beta}}\rho V^{\lambda}U^{\lambda}\wedge 
d\bar{u}_{m+1}\wedge \cdots \wedge d\bar{u}_p \wedge \hat{\varphi},
\end{equation}
where we have introduced the notations: $\alpha=\sum_1^p\alpha(j)$, $\beta=\sum_1^q\beta(j)$,
$\lambda\gamma=\sum_{m+1}^p\lambda_j\alpha(j)+\sum_1^q\lambda_{p+j}\beta(j)$, and
\begin{equation*}
V^{\lambda}=\frac{|v_1|^{2\lambda_{p+1}}\cdots |v_q|^{2\lambda_{p+q}}}{(v_1\cdots v_q)^N},
\end{equation*}
\begin{equation*}
U^{\lambda}=\lambda_{m+1}\cdots \lambda_p
\frac{|u_{m+1}|^{2\lambda_{m+1}-2}\cdots |u_{p}|^{2\lambda_{p}-2}}{
(u_{m+1}\cdots u_p)^{N-1}}.
\end{equation*}
Let $K\subseteq \{1,\ldots,n\}$ be the set of indices $i$ such that $x_i$ divides at least 
some $\hat{g}_j$. 
We will use the following division lemma, proved below, to replace the form 
$d\bar{u}_{m+1}\wedge \cdots \wedge d\bar{u}_p \wedge \hat{\varphi}$ in \eqref{prop1}
by another one, which vanishes on the variety $\{\prod_{i\in K}x_i=0\}$.

\begin{lemma}\label{divlemma}
If $\psi$ is a holomorphic $n-p$-form on the base manifold $X$, then one can find  
explicitly a holomorphic $n-m$-form $\omega$ in the $x$-chart on $\mathcal{X}$ such that
\begin{itemize}
\item[i)] $\frac{dx_j}{x_j}\wedge (du_{m+1}\wedge\cdots\wedge du_p\wedge\hat{\psi}-\omega)$ 
is non-singular for all $j\in K$, and
\item[ii)] $dx^{\alpha(1)}\wedge\cdots\wedge dx^{\alpha(m)}\wedge \omega =0$.
\end{itemize}  
\end{lemma}

\noindent By linearity, we may assume that $\varphi$ is decomposable and write
$\varphi=\phi_1\wedge \bar{\phi}_2$, where $\phi_1\in \D_{n,0}(X)$ and $\phi_2$ is a holomorphic
$n-p$-form. With $\phi_2$ as in-data to Lemma \ref{divlemma} we thus see that we may replace
$d\bar{u}_{m+1}\wedge \cdots \wedge d\bar{u}_p \wedge \hat{\varphi}$ in \eqref{prop1} by a 
$(n,n-m)$-form $\xi$, without affecting the integral, such that 
$(d\bar{x}_j/\bar{x}_j)\wedge \xi$ is smooth for all $j\in K$. It follows that for any 
$L\subseteq K$, $\wedge_{j\in L}(d\bar{x}_j/\bar{x}_j)\wedge \xi$ is a smooth form. 
Using Leibniz' rule to expand the expressions $\debar |x^{\alpha(j)}|^{2\lambda_j}$,
$1\leq j \leq m$, the integral 
\eqref{prop1} can be written
\begin{equation}\label{prop2}
\sum_{i_1<\cdots <i_m}
\det (A(i_1,\ldots,i_m))
\int \frac{|x^{\lambda(\alpha+\beta)}|^2}{x^{N(\alpha+\beta)}}
\frac{d\bar{x}_{i_1}\wedge \cdots \wedge d\bar{x}_{i_m}}{\bar{x}_{i_1}\cdots \bar{x}_{i_m}}
\wedge \Phi(x;\lambda),
\end{equation}
where $A(i_1,\ldots,i_m)$ is the matrix $\big( \alpha(i_k)_{i_l}\big)_{k,l}$, 
$\lambda(\alpha+\beta)=\sum_1^p\lambda_j\alpha(j)+\sum_1^q\lambda_{p+j}\beta(j)$, and 
\begin{equation*}
\Phi(x;\lambda)=
\lambda_1\cdots \lambda_m \rho V^{\lambda} U^{\lambda}\wedge \xi.
\end{equation*}
We emphasize that $\Phi$ is a smooth compactly supported form depending analytically on $\lambda$
and that $\wedge_{j\in L}(d\bar{x}_j/\bar{x}_j)\wedge \Phi$, $L\subseteq K$, also is. For notational
convenience, we consider the term of \eqref{prop2} with $i_j=j$ and to make our considerations 
non-trivial we then assume that $A:=A(1,\ldots,m)$ is non-singular. Furthermore, we assume, also
for simplicity, that 
$1,\ldots,k\notin K$ and that $k+1,\ldots,m\in K$. If we put $\mu=\lambda(\alpha+\beta)$ we can
write the term under consideration as
\begin{equation}\label{prop3}
\frac{\det (A)}{\mu_1\cdots \mu_k}
\int \frac{\debar |x_1|^{2\mu_1}\wedge \cdots \wedge \debar |x_k|^{2\mu_k}
|x_{k+1}|^{2\mu_{k+1}}\cdots |x_n|^{2\mu_n}}{x^{N(\alpha+\beta)}}\wedge
\end{equation}
\begin{equation*}
\hspace{5cm} \wedge
\frac{d\bar{x}_{k+1}\wedge \cdots \wedge d\bar{x}_{m}}{\bar{x}_{k+1}\cdots \bar{x}_{m}}
\wedge \Phi(x;\lambda).
\end{equation*}
Here, the expression on the second row is a smooth compactly supported form depending analytically 
on $\lambda$. After this observation it is a one-variable problem to see that 
the integral, without the coefficient in front, has an analytic continuation to some half space $H$ 
independent of $N$; see, e.g., Lemma 2.1 in \cite{matsa}. The possible poles are therefore only along
hyperplanes of the form $\mu_j=0$. Fix a $j$ with $1\leq j \leq k$. Then $j\notin K$, which 
means that $x_j$ does not divide any of the $\hat{g}$-functions. Hence,
$\beta(1)_j=\cdots =\beta(q)_j=0$, and consequently, $\mu_j=\sum_{i=1}^p\alpha(i)_j\lambda_i$.
Moreover, if $\mu_j$ happens to be proportional to some $\lambda_i$ then, first of all, 
$x_j$ must divide $x^{\alpha(i)}$ but no other $x^{\alpha(l)}$ (or any $x^{\beta(l)}$). 
Secondly, the term \eqref{prop3} of \eqref{prop2} that we are considering must have arisen from 
the term in the Leibniz expansion of \eqref{prop1} when the $\debar$ in front of 
$|x^{\alpha(i)}|^{2\lambda_i}$ has fallen on $|x_j^{\alpha(i)_j}|^{2\lambda_i}$. 
Thus, $i\leq m$ and no other $\mu_{\nu}$ with $1\leq \nu \leq k$ can be proportional 
to $\lambda_i$. Since $\Phi(x;\lambda)$ is divisible by $\lambda_i$ we can therefore cancel 
poles along hyperplanes $\mu_j=0$ if $\mu_j$ is proportional to some $\lambda_i$.
In conclusion, $\Gamma_{f,g}^{\varphi}(\lambda)$ has a meromorphic continuation to some 
half space $H$ with possible poles only along hyperplanes of the form
$\sum_1^pa_j\lambda_j=0$, where $a_j\in \N$ and at least two $a_j$ must be non-zero. 
\end{proof}

\begin{proof}[Proof of Lemma \ref{divlemma}]
Put $\Psi=du_{m+1}\wedge\cdots\wedge du_p\wedge\hat{\psi}$ and define
\begin{equation*}
\omega=\sum_{j\in K}\Psi_{j}-\sum_{\stackrel{i,j\in K}{i<j}}\Psi_{ij}+\cdots
+(-1)^{|K|-1}\Psi_{i_1\cdots i_{|K|}},
\end{equation*}
where $\Psi_{i_1\cdots i_{\ell}}$ means the pullback of $\Psi$ to $\{x_{i_1}=\cdots=x_{i_{\ell}}=0\}$
extended constantly to $\C^n$. A straight forward induction over $|K|$ shows that $\omega$
so defined satisfies i). (See also \cite{HS}.) To see that $\omega$ satisfies ii), consider 
a $\Psi_{i_1\cdots i_{\ell}}$. Let $L$ be the set of indices $j$ such that no $x_{i_{k}}$,
$1\leq k \leq \ell$, divides $\hat{f}_j$ and write $L=L'\cup L''$, where $L'=\{j\in L;\, j\leq m\}$
and $L''=\{j\in L;\, m+1\leq j\leq p\}$. For each $x_{i_k}$, with $1\leq k \leq \ell$,
we know that $x_{i_k}$ divides some $\hat{g}$-function, say $\hat{g}_{j_k}$. The variety
$\{x_{i_1}=\cdots = x_{i_{\ell}}=0\}$ is then contained in 
$\{\hat{g}_{j_1}=\cdots = \hat{g}_{j_{\ell}}=0\}\bigcap \cap_{i\notin L}\{\hat{f}_i=0\}$, i.e., in 
the preimage of $V:=\{g_{j_1}=\cdots = g_{j_{\ell}}=0\}\bigcap \cap_{i\notin L}\{f_i=0\}$. 
Since $(f,g_j)$ defines a complete 
intersection for any $j$, the variety $V$ has codimension at least $p-|L|+1$. 
Now, the form $\wedge_{j\in L}df_j\wedge \psi$ has degree $n-p+|L|$ and thus, has a vanishing 
pullback to $V$. Hence, we get that
\begin{equation*}
\widehat{\bigwedge_{j\in L}df_j\wedge \psi}=
\bigwedge_{j\in L}d\hat{f}_j \wedge \hat{\psi}=
\bigwedge_{i\in L'}dx^{\alpha(i)} \wedge \bigwedge_{j\in L''}d(u_jx^{\alpha(j)}) \wedge \hat{\psi}
\end{equation*}
has a vanishing pullback to $\{x_{i_1}=\cdots=x_{i_{\ell}}=0\}$.
But this means that 
\begin{equation*}
x^{\sum_{i\in L''}\alpha(i)}\bigwedge_{j\in L'} dx^{\alpha(j)}\wedge
(\bigwedge_{k\in L''}du_{k}\wedge \hat{\psi})_{i_1\cdots i_{\ell}}+
\bigwedge_{\iota\in L'} dx^{\alpha(\iota)}\wedge
\sum_{\nu\in L''}dx^{\alpha(\nu)}\wedge \xi_{\nu}=0,
\end{equation*}
for some forms $\xi_{\nu}$, where the first term arises when no differential 
hits any $x^{\alpha(j)}$, $j\in L''$. Taking 
the exterior product with $\wedge_{j\notin L''}(du_j)_{i_1\cdots i_{\ell}}$ we obtain
\begin{equation*}
x^{\sum_{i\in L''}\alpha(i)}\bigwedge_{j\in L'} dx^{\alpha(j)}\wedge
\Psi_{i_1\cdots i_{\ell}}+\bigwedge_{\iota\in L'} dx^{\alpha(\iota)}\wedge
\sum_{\nu\in L''}dx^{\alpha(\nu)}
\wedge\tilde{\xi}_{\nu}=0.
\end{equation*}
We now multiply this equation with the exterior product of all $dx^{\alpha(j)}$ with 
$j\leq m$ and $j\notin L'$. Then we get $dx^{\alpha(1)}\wedge \cdots \wedge dx^{\alpha(m)}$ in front
of the sum and this makes all terms under the summation sign disappear since every $\alpha(\nu)$,
with $\nu \in L''$, is a linear combination of $\alpha(1),\ldots,\alpha(m)$. It thus follows that 
\begin{equation*}
x^{\sum_{i\in L''}}dx^{\alpha(1)}\wedge \cdots \wedge dx^{\alpha(m)}
\wedge \Psi_{i_1\cdots i_{\ell}}=0,
\end{equation*}
and since this holds everywhere we may remove the factor $x^{\sum_{i\in L''}}$ and conclude that 
$\omega$ has the property ii).
\end{proof}

\begin{proof}[Proof of Theorem \ref{analytisk}]
The proof is based on induction over $p$. The induction start, $p=1$, follows from Proposition
\ref{polprop}. Assume therefore that the theorem is proved for $p=k$. We introduce the notation
$\gamma(\lambda_{i_1},\ldots,\lambda_{i_p};\lambda_*)$
for the current-valued function
\begin{equation*}
\frac{\debar |f_{i_1}|^{2\lambda_{i_1}}\wedge \cdots \wedge \debar |f_{i_p}|^{2\lambda_{i_p}}
|f_{j_1}|^{2\lambda_{j_1}}\cdots |f_{j_q}|^{2\lambda_{j_q}}}{f_1^N\cdots f_{p+q}^N}.
\end{equation*}
When all $\real \, \lambda_j$ are large we have the equality of currents
\begin{equation*}
\debar \gamma(\lambda_1,\ldots,\lambda_k;\lambda_*) =
(-1)^k\sum_{j=1}^q\gamma(\lambda_1,\ldots,\lambda_k,\lambda_{k+j};\lambda_*),
\end{equation*}
and by our assumption, the left hand side is analytic in some half space $H$ independent of $N$.
Moreover, by Proposition \ref{polprop}, the term on the right hand side corresponding to $j$
has only poles along hyperplanes of the form $a_{k+j}\lambda_{k+j}+\sum_1^ka_i\lambda_i=0$.
It thus follows that $\gamma(\lambda_1,\ldots,\lambda_{k+1};\lambda_*)$, on one hand, 
only has poles along hyperplanes $\sum_1^{k+1}a_i\lambda_i=0$ and, on the other,
only has poles along hyperplanes $a_{k+j}\lambda_{k+j}+\sum_1^{k}a_i\lambda_i=0$ with $j>1$. 
But then $\gamma(\lambda_1,\ldots,\lambda_{k+1};\lambda_*)$ can only have poles along hyperplanes 
of the form $\sum_1^{k}a_i\lambda_i=0$. Consider now the current equality
\begin{eqnarray*}
\debar \gamma(\lambda_1,\ldots,\lambda_{k-1},\lambda_{k+1};\lambda_*) &=&
(-1)^{k+1}\gamma(\lambda_1,\ldots,\lambda_{k+1};\lambda_*) \\
&+&
(-1)^k\sum_{2}^q\gamma(\lambda_1,\ldots,\lambda_{k-1},\lambda_{k+1},\lambda_{k+i};\lambda_*).
\end{eqnarray*}
From this it follows similarly that $\gamma(\lambda_1,\ldots,\lambda_{k+1};\lambda_*)$ can only 
have poles along hyperplanes of the form $a_{k+1}\lambda_{k+1}+a_{k+j}\lambda_{k+j}+
\sum_1^{k-1}a_i\lambda_i=0$ with $j> 1$. Since we know that its only poles are along 
$\sum_1^{k}a_i\lambda_i=0$, we see that it in fact only can have poles along
$\sum_1^{k-1}a_i\lambda_i=0$. Continuing in this way, looking at appropriate current 
equalities and using the induction hypothesis and Proposition \ref{polprop},
we eventually see that $\gamma(\lambda_1,\ldots,\lambda_{k+1};\lambda_*)$ cannot have any poles 
at all in $H$. This concludes the induction step and consequently the proof of Theorem \ref{analytisk}.
\end{proof}

\end{document}